**Словак К. І., Семеріков С. О., Триус Ю. В.**
Криворізький національний університет,
Черкаський державний технологічний університет

# Мобільні математичні середовища: сучасний стан та перспективи розвитку

**Вступ.** На формування і розвиток особистості людини найбільше впливає середовище, в якому вона живе, навчається, працює. Тому для ВНЗ важливою і актуальною проблемою як теоретичного, так і практичного характеру є проблема створення такого високотехнологічного інформаційно-комунікаційного освітньо-наукового середовища, в якому студент знаходиться щодня в процесі всього періоду навчання у вищій школі, яке повинне відповідати потребам інформаційного суспільства, сучасному рівню науки, техніки та світовим освітнім стандартам і сприяти підвищенню рівня інформаційно-комунікаційної підготовки та формуванню професійних компетентностей. Прикладом такого інформаційно-комунікаційного освітньо-наукового середовища для навчання математичних дисциплін студентів ВНЗ є web-орієнтоване математичне середовище. Одному з підходів до створення такого середовища на основі web-орієнтованої системи комп'ютерної математики (СКМ), що інтегрує в собі послуги різних систем за допомогою клієнт-серверних технологій і таких засобів ІКТ навчання математики, як мультимедійні демонстрації, динамічні математичні моделі, тренажери та експертні системи навчального призначення, присвячена ця стаття. Зокрема в роботі визначено новий клас педагогічного програмного забезпечення математичних дисциплін – *мобільні математичні середовища*.

Однією з проблем, що постає при створенні математичного середовища, є вибір системи комп'ютерної математики, що буде ядром цього середовища. Як комерційні, так і вільно поширювані системи комп'ютерної математики суттєво відрізняються за функціональністю (загального призначення, спеціалізовані), набором допустимих математичних і графічних операцій, засобами програмування і управління даними, операційною платформою, швидкістю обчислень, інтерфейсом (командного рядка, графічним, Web-), розміром (від кількасот кілобайт до кількох гігабайт), зручністю у навчанні і використанні тощо. Крім того, суттєвою вимогою до ядра математичного середовища є можливість інтегрувати його з іншими системами комп'ютерної математики, що обумовлено, насамперед, необхідністю:

– раціонального вибору СКМ з урахуванням особливостей задачі, що розв'язується;
– розв'язування складних задач за допомогою різних засобів, щоб перевірити правильність результатів, не покладаючись на одну СКМ (збільшити вірогідність одержаного результату);
– підготовки математичних документів навчального призначення.
– Спільність функціональних характеристик різних СКМ може спонукати до розробки деякої «універсальної» СКМ за одним з напрямів:
– створення мінімального ядра, спільного для всіх систем («переріз»);
– створення максимального ядра, до якого входять засоби всіх систем («об'єднання»);
– створення мінімального ядра, до якого за необхідністю інтегруються інші системи, визначені розробником чи користувачем.

За першим та другим напрямом створюються нові СКМ – або занадто обмежені, або надмірно громіздкі. Останнє говорить на користь третього напряму – інтеграції математичних систем між собою і з іншими програмами, що може розглядатися як один з перспективних напрямів розвитку систем комп'ютерної математики [1, 364–365].

Розглянемо детальніше основні характеристики і структуру *мобільного математичного середовища (ММС),* критерії вибору СКМ для ядра ММС, проблеми створення навчально-методичної складової інформаційного забезпечення ММС, а також приклад реальної ММС «Вища математика».

**Основна частина**. *Мобільне математичне середовище* – відкрите модульне мережне мобільне інформаційно-обчислювальне програмне забезпечення, що надає користувачеві (викладачеві, студентові) можливість мобільного доступу до інформаційних ресурсів математичного і навчального призначення, умови для ефективної організації процесу навчання та інтеграції аудиторної і позааудиторної роботи.

До основних характеристик ММС відносяться [2; 3]:

– *мобільність доступу*: доступ з широкого спектру комп'ютерних пристроїв, що надає можливість залучити в якості засобів навчання нетбуки, планшетні комп'ютери та смартфони;
– *мобільність програмного забезпечення*: можливість перенесення середовища на різні програмно-апаратні платформи без суттєвої модифікації;

- *мережність*: зберігання математичних об'єктів на мережних серверах, що надає можливість уніфікувати доступ до них як в навчальній аудиторії, так і за її межами;
- *відкритість*: можливість зміни інформаційної та обчислювальної складової середовища;
- *модульність*: можливість додавання, вилучення та заміни компонентів середовища;
- *об'єктна орієнтованість*: можливість прототипування, створення, модифікації, наслідування, інкапсуляції математичних об'єктів;
- *можливість* природного застосування ефективних педагогічних технологій *організації спільної роботи* над навчальними проектами у навчальних спільнотах.

Особливістю ММС є динамічна природа навчальних матеріалів – будь-який опублікований у мережі об'єкт може автоматично змінюватися у відповідності до: зміни вмісту пов'язаного з ним робочого аркуша; зміни програмного забезпечення, що входить до складу ММС; зміни пристрою доступу до навчальних матеріалів; зміни початкових умов для моделей.

*Основними складовими* ММС є *обчислювальне ядро* (математичний пакет – інтегратор СКМ) та *інформаційне забезпечення*, що містить навчально-методичні та допоміжні інформаційні матеріали.

У процесі створення ММС для навчання математичних дисциплін особливу увагу слід приділити вибору математичного пакету, що складає ядро ММС.

Головними критеріями вибору СКМ для ядра ММС є:
- розширюваність (система повинна легко модифікуватися користувачем за рахунок її доповнення для розв'язування нових класів задач);
- наявність різних інтерфейсів та підтримка web-сервісів (для забезпечення мобільного доступу);
- кросплатформеність (мобільність програмного забезпечення);
- можливість створення програм із стандартними елементами управління (лекційних демонстрацій, динамічних моделей, тренажерів, навчальних експертних систем);
- можливість інтегрувати в ММС різноманітне програмне забезпечення (на основі відкритих програмних інтерфейсів);
- підтримка технології Wiki;
- можливість локалізації та вільне поширення.

Як показали дослідження авторів [2-6], найбільш повно наведеним критеріям відповідає Web-СКМ Sage (таблиця 1), використання якої надає можливість в рамках одного середовища реалізувати основні складові ММС, зокрема, лекційні демонстрації, динамічні моделі, тренажери, навчальні експертні системи та інші навчальні матеріали; автоматизувати обчислювальний процес розв'язування задач прикладної спрямованості, зосередившись на побудові математичної моделі та інтерпретації результатів обчислювального експерименту; інтегрувати аудиторну та позааудиторну роботу студентів в єдиному середовищі. Саме тому в якості обчислювального ядра ММС доцільно обрати Web-СКМ Sage.

Враховуючи, що інформаційне забезпечення, яке входить до складу ММС, є предметно-орієнтованим, розглядатимемо клас мобільних математичних середовищ, що мають однакове обчислювальне ядро та варіативне інформаційне забезпечення (рис. 1). Таким чином, заміна методичної складової інформаційного забезпечення надає можливість створювати різні ММС навчання математичних дисциплін.

Розглянемо ММС «Вища математика», призначене для підтримки навчання вищої математики студентів економічних спеціальностей [7].

Розробка методичної складової інформаційного забезпечення ММС «Вища математика» здійснювалась за такими напрямами: 1) графічна інтерпретація математичних моделей та теоретичних понять; 2) автоматизація рутинних обчислень; 3) підтримка самостійної роботи; 4) математичні дослідження; 5) генерація навчальних завдань. При цьому перші чотири напрями спрямовані на активізацію навчальної діяльності студентів, а п'ятий – на підвищення ефективності діяльності викладача.

Для реалізації *першого* та *четвертого* напрямів створено комп'ютерні моделі з графічним інтерфейсом і напівавтоматичним режимом управління. Розроблені моделі різняться за дидактичним призначенням відповідно до вказаних напрямів: лекційні демонстрації – для унаочнення абстрактних математичних понять, динамічні моделі – для проведення навчальних досліджень. В розроблених моделях можуть використовуватися Web-сервіси доступу до баз знань (наприклад, Wolfram|Alpha) та баз даних (зокрема, Google Finance).

Реалізація *другого* напряму передбачала використання, насамперед, обчислювального ядра ММС з метою автоматизації процесу дослідження задач з економічним змістом, зосереджуючись на побудові моделі та інтерпретації результатів обчислювального експерименту.

Для реалізації *третього* напряму розроблено:

– індивідуальні домашні завдання до кожного модуля, що включають приклади розв'язування типових завдань за темою модуля та задач для самостійного опрацювання трьох рівнів (для відпрацювання навичок «ручного» розв'язування; комп'ютерно-орієнтовані задачі, витрати часу на «ручне» розв'язування яких перевищують час на створення та дослідження моделі; дослідницькі завдання);

– приклади розв'язування завдань з кожного модуля у традиційному вигляді та за допомогою ММС (при цьому особливості компонування завдань, детальні пояснення кожного кроку розв'язування, застосування засобів ІКТ сприяють підвищенню ефективності самостійної роботи студентів);

– програми-тренажери з покроковою деталізацією етапів розв'язування математичної задачі, що надає можливість студентам здійснити детальну перевірку кожного кроку виконання завдання;

– навчальні експертні системи, використання яких надає викладачеві можливість організувати автоматизований контроль та корекцію результатів навчальної діяльності студентів, проводити тренування та підготовку до модульного і підсумкового контролю. Самостійна робота студентів зі створення баз знань експертної системи навчального призначення за обраною темою курсу вищої математики сприяє узагальненню та систематизації знань;

– навчально-консультативний форум надає можливість для спілкування всім користувачам різних інсталяцій ММС «Вища математика», дистанційного консультування, обміну досвідом розв'язування задач засобами ММС.

*Таблиця 1*

Порівняльна характеристика СКМ

| Математичний пакет | Ліцензія | Розширюваність | Web-інтерфейс | WAP-інтерфейс | XML-RPC/SOAP/інші Web-сервіси | Кросплатформеність | Можливість створення програм із стандартними елементами управління | Можливість інтеграції різних пакетів | Підтримка Wiki | Інтернаціоналізація та візуалізація | Мови програмування |
|---|---|---|---|---|---|---|---|---|---|---|---|
| GRAN | ВП, ЗК | | | | | | + | | | + | |
| DG | ВП, ЗК | | | | | | + | | | | |
| GeoGebra | ВП, ВК | + | + | | | + | | | | + | GGBScript, JavaScript |
| KAlgebra | ВП, ВК | + | | | | + | | | | + | |
| Cinderella | КП, ЗК | + | + | | | + | + | | | + | CindyScript |
| Maxima | ВП, ВК | + | + | + | | + | | | | + | Maxima, Lisp |
| Sage | ВП, ВК | + | + | + | + | + | + | + | + | + | Python, Sage, Sh, GAP, GP, Maxima, R, Singular, Kash, Macaulay, Magma, Maple, Mathematica, MATLAB, Mupad, Octave, Scilab, Spad |
| Maple | КП, ЗК | + | + | | + | + | + | + | | | Maple, C/C++, Fortran, MATLAB |
| MATLAB | КП, ЗК | + | + | | + | + | + | + | | | MATLAB, C/C++, Fortran, Java |
| Mathematica | КП, ЗК | + | + | | + | + | + | + | | | Mathematica, C/C++, Fortran, Assembler |
| Magma | КП, ЗК | + | + | | | + | | | | | Magma |
| MathCAD | КП, ЗК | + | + | | + | + | + | + | | + | |

*Позначення:*

ВП – вільне поширення, КП – комерційне поширення, ВК – відкриті коди, ЗК – закриті коди

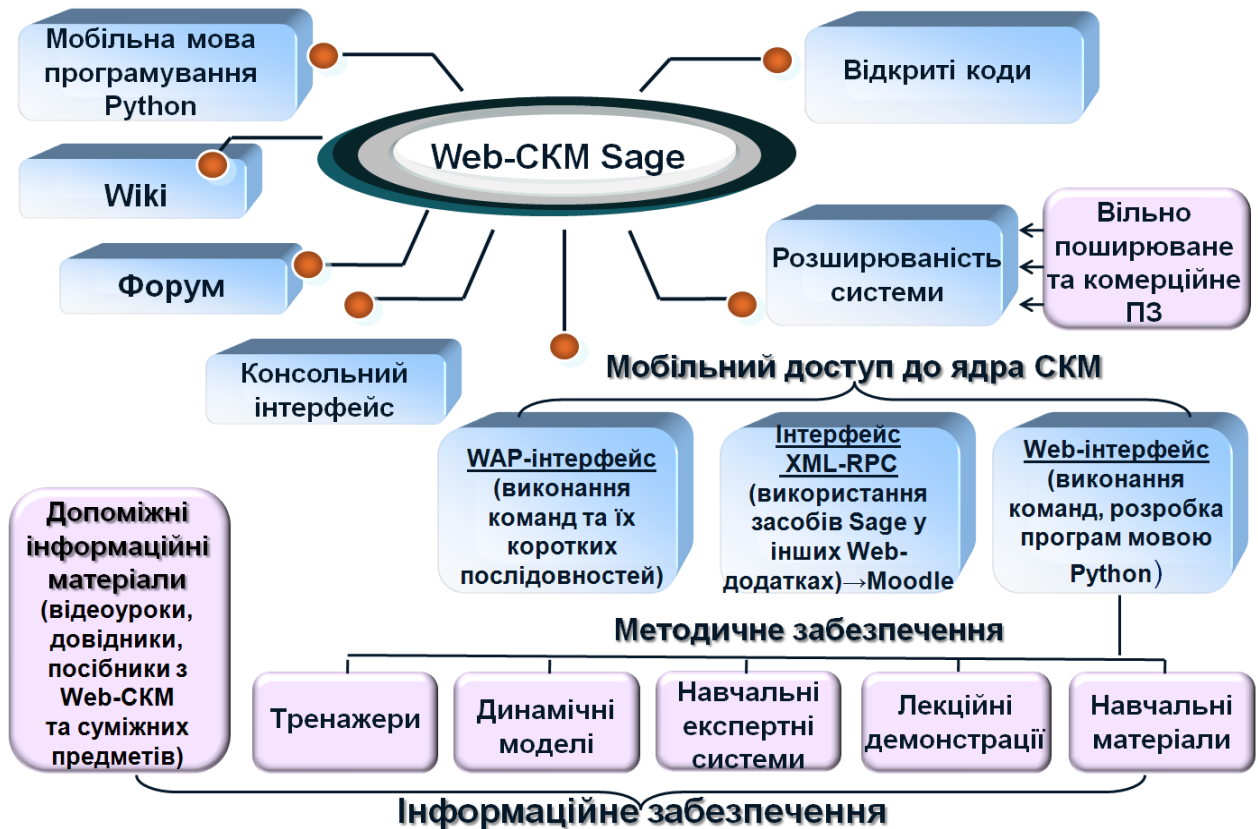

Рис. 1. Загальна структура мобільного математичного середовища

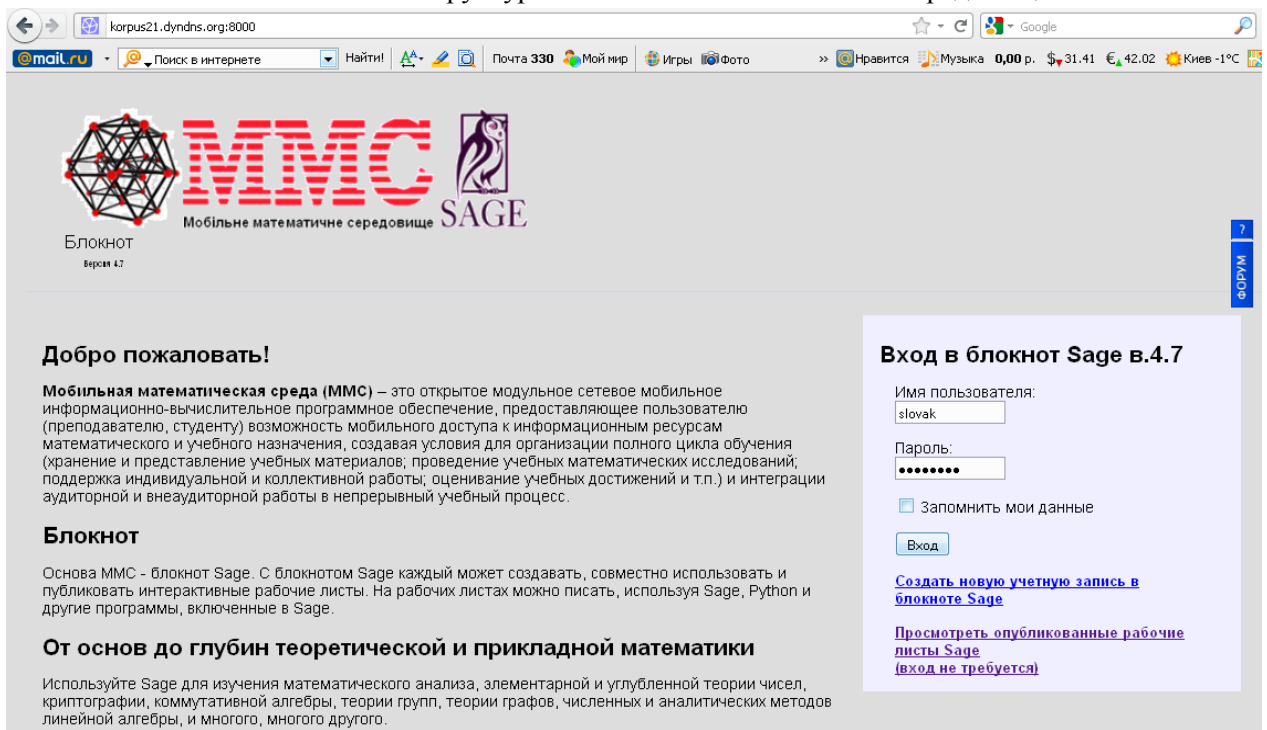

Рис. 2. Web-сервер ММС «Вища математика»

Для автоматизації підготовки та перевірки індивідуальних домашніх завдань, матеріалів для поточного та модульного контролю створено програми-генератори навчальних завдань, для виконання завдань *п'ятого* напряму. Завдяки можливості збереження результату генерації у природній математичній нотації, вибору довільної кількості завдань із відповідями, ММС є ефективним засобом створення генераторів навчальних завдань.

Для забезпечення відкритого доступу до розробленого ММС «Вища математика» було налаштовано Web-сервер ММС за адресою http://korpus21.dyndns.org:8000/ (рис. 2). Кожному зареєстрованому користувачеві автоматично завантажуються всі навчально-методичні матеріали ММС «Вища математика» (рис. 3), крім індивідуальних домашніх завдань, які потрібно обрати самостійно із списку опублікованих аркушів згідно номеру модуля, що вивчається, та номеру студента в журналі академічної групи (рис. 4).

Назви робочих аркушів ММС «Вища математика» подано у форматі: `M_номер_призначеннння_назва (номер)`, де

M – «модуль»;

`номер` – означає, до якого саме змістового модуля курсу вищої математики належить робочий аркуш (набуває значень від 1 до 9);

`призначення` – вказує на тип програмного модуля чи розділ методичного комплексу (D – лекційні демонстрації та динамічні моделі, G – генератори навчальних завдань, T – тренажери, ES – навчальні експертні системи, Ex – вправи, ExS – приклади розв'язування, L – лекції, ІДЗ_варіант – варіант індивідуального домашнього завдання (від 1 до 30)).

Для підвищення зручності використання ММС «Вища математика» до його складу включено робочі аркуші з основними відомостями щодо організації роботи у Web-СКМ Sage (рис. 3). Аркуші містять теоретичні відомості, приклади застосування та відеоуроки. Крім того, з метою усунення незручностей, пов'язаних з якістю послуг, що надаються провайдерами Інтернет, кожен користувач може завантажити локальну Windows-версію ММС «Вища математика», призначену для роботи в автономному режимі. Завдяки цьому швидкість виконання обчислень у ММС, пропорційна потужності комп'ютерної системи та обсягу оперативної пам'яті, може бути вище, ніж на сервері ММС, що економить навчальний час (це особливо актуально для лекційних занять). За наявності з'єднання з мережею Інтернет з локальної версії ММС можна звернутися до форуму. Для завантаження локальної Windows-версії ММС «Вища математика» необхідно подвійним натисканням правої кнопки миші відкрити закладку «Форум» та знайти на ньому відповідні посилання (рис. 5). Зауважимо, що використання лише локальної версії ММС «Вища математика» обмежує можливості у спілкуванні студентів з викладачем та суттєво ускладнює оновлення навчально-методичних матеріалів.

Рис. 3. Навчально-методичні матеріали ММС «Вища математика»

Рис. 4. Список опублікованих ІДЗ ММС «Вища математика»

Рис. 5. Навчально-консультативний форум ММС «Вища математика»

**Висновки**. 1. Мобільне математичне середовище є інноваційним засобом навчання математичних дисциплін, визначальними особливостями якого є: об'єднання інших засобів навчання (лекційних демонстрацій, динамічних моделей, тренажерів та навчальних експертних систем); можливість налаштування на конкретну математичну дисципліну; динамічна природа навчальних матеріалів.

2. Перспективними напрямами розвитку ММС є: розробка засобів доступу до ММС з мобільних пристроїв на платформах Google Android та Apple IOS; вбудовування робочих аркушів ММС та їх компонентів у системи підтримки навчання; інтеграція з соціальними мережами; розробка засобів доступу до Wolfram|Alpha.

3. Перспективами розвитку методики використання ММС є: розробка ММС з інших математичних дисциплін для студентів ВНЗ з урахуванням професійної спрямованості навчання; розробка методики використання засобів ММС у процесі навчально-дослідницької роботи студентів: підготовка конкурсних, курсових, кваліфікаційних робіт молодшого спеціаліста, бакалавра, спеціаліста, магістра; розробка методики використання засобів ММС у процесі навчання математичних спецкурсів за різними напрямами підготовки.

**Література**